\documentclass[11pt]{article}
\usepackage{amsmath}
\usepackage{amsthm}
\usepackage{graphicx}
\usepackage{bbm,amssymb,enumerate}
\usepackage[hyperfootnotes=false]{hyperref}
\usepackage{float}

\def\R{\mathbb R}
\def\Z{\mathbb Z}

\def\O{\Omega}
\def\disp{\displaystyle}

\newcommand{\old}[1]{}
\newcommand{\arxiv}[1]{{\tt \href{http://arxiv.org/abs/#1}{arXiv:#1}}}

\newtheorem{theorem}{Theorem}

\newtheorem{lemma}[theorem]{Lemma}

\theoremstyle{remark}

\theoremstyle{definition}

\def \a {\alpha}

\def\g{\gamma}

\def \d{\delta}
\def \D{\Delta}
\def\e{\epsilon}
\def\f{\varphi}
\def\F{\Phi}

\def \l {\lambda}

\def\o{\omega}
\def \O {\Omega}
\def\p{\psi}

\def \s {\sigma}

\def\ll{\mathcal L}
\def\ee{\mathbf E}
\def\pp{\mathbf P}
\def\Z{\mathbf Z}
\def\R{\mathbf R}

\def\T{\mathbf T}
\def\back{\backslash}
\def\disp{\displaystyle}

\title{Internal DLA for Cylinders}
\author{David Jerison\footnote{Partially supported by NSF grant DMS-1069225.} 
\and Lionel Levine\footnote{Partially supported by NSF grant DMS-1105960.} \and Scott Sheffield\footnote{Partially supported by NSF grant
DMS-0645585.}}
\date{May 30, 2012}

\begin{document}

\maketitle


\bigskip 

\begin{center}
\emph{Dedicated to E. M. Stein}
\end{center}

\bigskip

\section{Introduction}

Internal Diffusion-Limited Aggregation (internal DLA)  is a random
lattice growth model.  Consider  
the two-dimensional lattice, $\Z\times \Z$.   In the case
of a single source at the origin, the random 
occupied set $A(T)$ of $T$ lattice sites is defined inductively as follows.  
Let $A(1)$ be the singleton set containing the origin.  
Given $A(T-1)$, start a random walk in $\Z\times \Z$ at the origin.  
Then 
\[
A(T) : = \{n\} \cup A(T-1)
\]
where $n\in \Z\times \Z$ is the first site reached by the random walk that is not in $A(T-1)$.  

In this paper, we will discuss the continuum limit of internal DLA,
which is governed by a deterministic fluid flow equation known as 
Hele-Shaw flow.  Our main focus will be on 
fluctuations.  In \cite{JLS11} we characterized the average fluctuations of
the model just described in terms of a close relative of the Gaussian Free Field, defined below.  In this article we will prove the analogous results
for the lattice cylinder.  In the case of the cylinder, the fluctuations are
described in terms of the Gaussian Free Field exactly.
We will also state without proof an almost sure bound on the maximum 
fluctuation in the case of the cylinder analogous to the case of the planar
lattice proved in \cite{JLS}.    The 
main tools used in the proofs are martingales.  As we shall see,
the martingale property in this context is the counterpart in probability theory
of well-known conservation laws for Hele-Shaw flow.

\begin{figure}
\begin{center}
\includegraphics[height=.3\textheight]{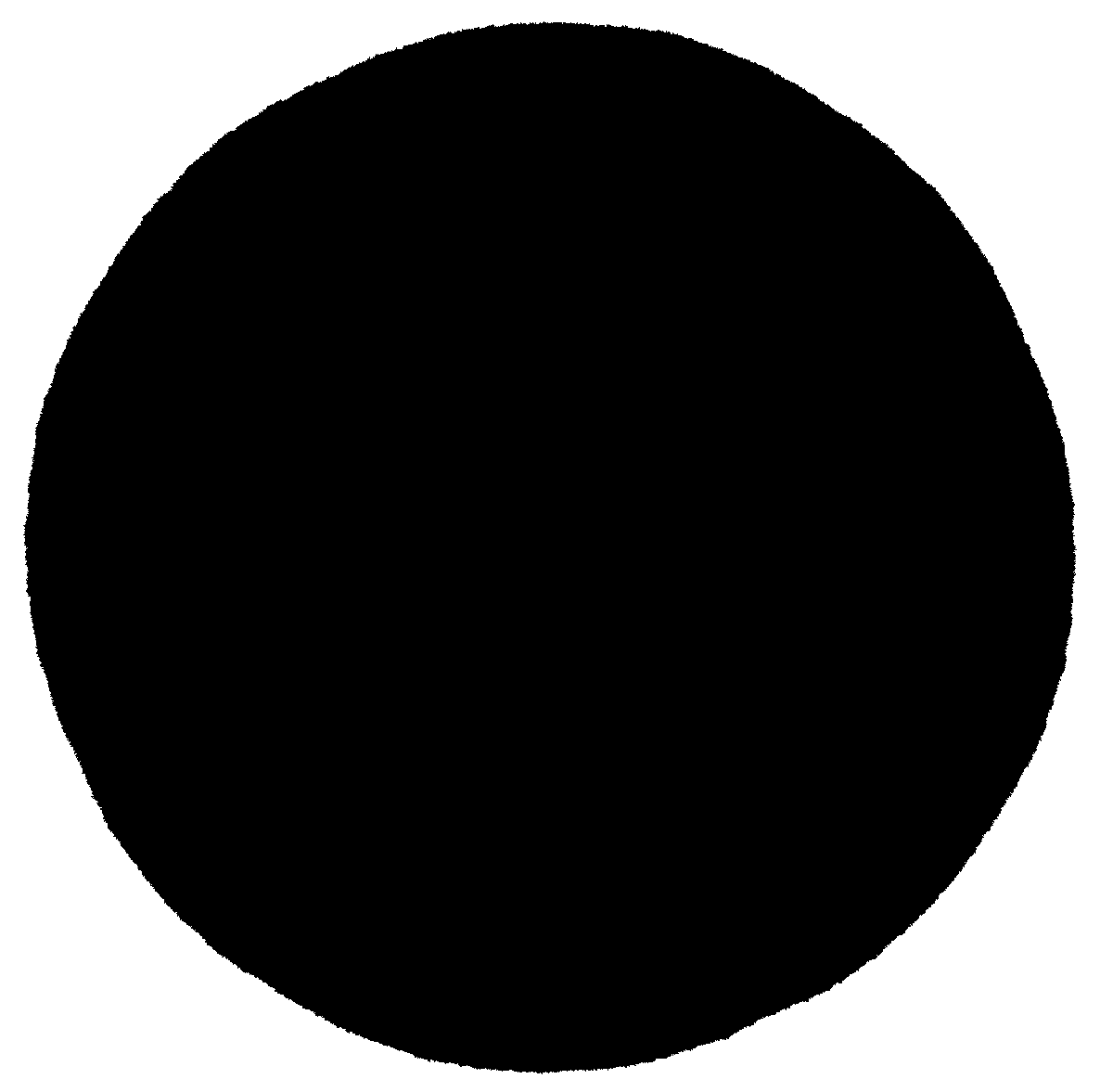} \\ \vspace{3mm}
\caption{\label{blob} 
Internal DLA cluster $A(T)$ with $T=10^6$ sites in $\Z^2$.}
\end{center}
\end{figure}

\begin{figure}
\begin{center}
\includegraphics[width=.7\textwidth]{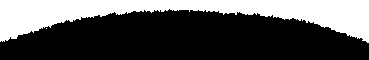} \\ \vspace{3mm}
\caption{\label{blob edge} 
Detail of boundary of the 1 million particle cluster or blob.
}
\end{center}
\end{figure}

The internal DLA model was introduced in 1986 by Meakin and Deutch \cite{MD} 
to describe chemical processes such as electropolishing, etching, and corrosion.  
Think of the occupied region as a blob of fluid.  Figure \ref{blob} depicts a simulation of a cluster (blob) of size one million in dimension $2$.  At each step a corrosive molecule is introduced at a source, which, in this simulation is a single point at the origin.
The corrosive particle wanders at random through the fluid until it
reaches the fluid-metal boundary, where it eats away a tiny portion of metal and enlarges slightly the fluid region.   The question that concerned Meakin and
Deutch was the smoothness of the surface that is being polished, that is,
how irregular the boundary is.  Figure \ref{blob edge} is a close-up picture of 
the boundary fluctuations.

Figure \ref{blob} suggests that the limit shape from a point source is a disk.    
Indeed, in 1992,
Lawler, Bramson and Griffeath \cite{LBG} proved that  the rescaled
limit shape of internal DLA from a point source is a 
ball in any dimension.  In 1995, Lawler \cite{Lawler95} proved
almost sure bounds on the cluster of the form 
\[
B(r - Cr^{1/3})\cap \Z^d \subset A(T) \subset B(r+Cr^{1/3}), \quad 
\]
where $T$ is the volume of the ball of radius $r$ and $C$ is a dimensional
constant.  
On the other hand, the numerical simulations of Meakin and Deutch predicted 
fluctuations, on average, of size $O(\sqrt{\log r})$ in dimension $2$ and
$O(1)$ in dimension $3$.   They made their predictions based on small values of $T$, but much larger simulations are now possible and give the same results.

The theorems we will describe are consistent 
with the size of fluctuations predicted by Meakin and Deutch  and reveal deeper structure, namely that the fluctuations
obey a central limit theorem.   The Fourier coefficients of the fluctuations tend to independent gaussians,
whose variance we can compute.   This gives a heuristic explanation of
numerical results on average fluctuations and many other predictions such as 
what should be the best possible bound on maximum fluctuations. 
In 2010, Asselah and Gaudilli\`ere \cite{AG10} improved the power in 
Lawler's bound in dimensions greater than $2$.   Later in 2010, Assellah and Gaudilli\`ere \cite{AG,AG2} and the present authors \cite{JLS,JLS2} independently proved logarithmic bounds on the maximum fluctuation. 

\begin{theorem} \label{maxfluc}  (Maximum Fluctuations)  There is a dimensional
constant $C_d$, such that  almost surely for sufficiently large $r$, 
\[
B(r- C_2 \log r) \cap \Z^2 \subset A(T) \subset B(r+ C_2 \log r)
\]
with $T = \pi r^2$.  Moreover, for $d\ge 3$, 
\[
B(r- C_d \sqrt{\log r}) \cap \Z^d \subset A(T) \subset B(r+ C_d \sqrt{\log r})
\]
where $T$ is the volume of the ball of radius $r$. 
\end{theorem}    
The maximum fluctuations represent the worst case along the entire circumference
as opposed to the average fluctuations observed by Meakin and Deutch.
Whether one considers the average or the worst case,
the model produces remarkably smooth surfaces --- 
even more smooth in dimension $3$ than in dimension $2$.
 
Before going any further, we should add a disclaimer.  Despite their superficial
similarity, the internal DLA model and the Diffusion-Limited Aggregation (DLA) model introduced by Witten and Sander \cite{WS} are very different.  
DLA is a model of particle deposition, in which a seed particle is placed at the 
origin in a lattice.  Particles follow a random walk starting at infinity and attach to the
existing cluster the first time they are adjacent to it.  The particles form a cluster
of fractal character and the continuum limit is very far from deterministic.
In their 1986 article, Meakin and Deutch refer to the work of Witten and Sander and explain that the internal DLA model is better behaved than DLA and intended
to describe quite different physical phenomena, ones that do not exhibit chaos.  
The Hele-Shaw model is also 
highly relevant to DLA, but it is the complement of the
cluster that is interpreted as the fluid region.  Thus the fluid region
shrinks.  When fluid is sucked away, the Hele-Shaw equation is
ill-posed, and the methods of partial differential equations no
longer apply except at very short time scales.  Instead, algebraic methods 
are used.   The subject is  
of great interest in statistical physics and has a direct connection with 
random matrices, but it is not the subject of this paper.

This paper discusses various aspects of several works of
the authors \cite{JLS,JLS2,JLS11}.  Rather than prove any
of the theorems in those papers, which concern $\Z^d$,
we prove two central limit theorems (Theorems \ref{circle theorem} and \ref{gff})
in which the set $\Z^2$ of \cite{JLS11} is replaced by 
the lattice cylinder $(\Z/N\Z)\times \Z$.   In the next section, we state 
our theorems in this new geometric setting.   In the third 
section we explain the 
relationship between internal DLA and Hele-Shaw flow.  
Sections 4 and 5 give complete proofs of two central limit
theorems for fluctuations of internal DLA on cylinders.  
We discuss the work of Levine and Peres concerning the
relationship of internal DLA with the obstacle problem in Section 6.
In the last section we make a few further remarks about the
theorems of \cite{JLS,JLS2,JLS11}, the 
effects of geometry on the problem, higher-dimensional questions,
and questions related to more general random walks.

\section{Main results for the cylinder.}

In this section we state our main results in the case of the two-dimensional cylinder
rather than the single source model in the plane which is carried out in \cite{JLS}.    We will make a comparison at the end of the paper.

Consider the cyclic group $\Z_N  = \Z/N\Z$, whose elements 
will typically be denoted  $n_1 = 1,2,\dots, N$.   In the
lattice cylinder $\Z_N\times \Z$, define the set
\[
A(0) = \{n=(n_1,n_2) \in \Z_N\times \Z: n_2 \le 0\}
\]
For integers $T>0$, the set $A(T)$ of lattice points is defined inductively,
with source at $n_2  = -\infty$.  Equivalently, 
given the set $A(T-1)$, start a random walk in $\Z_N\times \Z$ 
at one of the sites $(n_1,0)$, $n_1=1,\dots, N$, with equal probability.
$A(T)\back A(T-1)$ consists of the site at which the random
walk exits $A(T-1)$ for the first time.  Denote
\[
A^+(T) = A(T) \back A(0)
\]

A theorem analogous to Theorem \ref{maxfluc}, stated in a slightly more
precise form, is 
\begin{theorem}\label{maxfluccyl}  Given $0 <  y_1$ and $a <\infty$, there is a constant
$C$ depending only on $y_1$, and $a$ such that  with probability $1- N^{-a}$,  
for all $y$, $0 \le y \le y_1$, 
\[
\{n: n_2 \le yN - C\log N\}  \cap (\Z_N\times Z)
\subset 
A(T) \subset \{n: n_2 \le yN + C\log N\}
\]
with $T = \lfloor yN^2 \rfloor$. 
\end{theorem}

\begin{figure}
\begin{center}
\includegraphics[width=.8\textwidth]{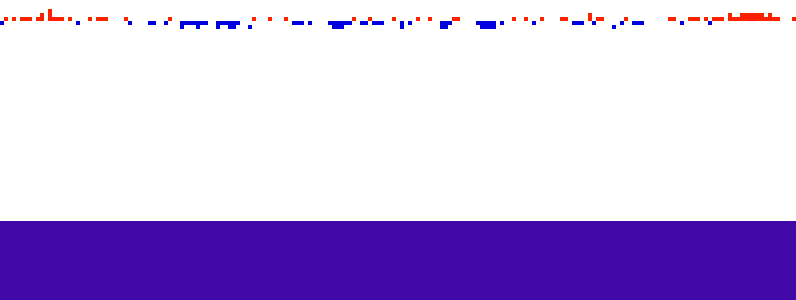} \\ \vspace{3mm}
\caption{\label{symm diff}
The symmetric difference of $A_N(T)$ and $\{y\le T/N^2 \}$ is the thin, ragged band
at the top with early points in red above the line $y=T/N^2$ and late points in blue below.  The bar at the bottom is the region $y\le 0$. 
}
\end{center}
\end{figure}

Next, we scale $A(T)$ by the factor $1/N$ to obtain a subset $A_N(T)$ 
of $\T \times \R$ with $\T = \R/\Z$.  
For $n= (n_1,n_2)$, $n_1 = 1,\dots, N$, $n_2\in \Z$, and
$0< x \le 1$ representing $x\in \T$, let
\begin{equation} \label{boxdef}
Q_N(n) = \{ (x,y)\in \T\times \R: n_1-1 <Nx \le n_1, \ n_2-1 <  Ny \le  n_2\}
\end{equation}
$Q_N(n)$ is the square of sidelength $1/N$ with $n/N$ at its upper right 
corner.    
Define \begin{equation} \label{AN(T)}
A_N(T) = \bigcup_{n\in A(T)} Q_N(n); \quad 
A_N^+(T) = \bigcup_{n\in A^+(T)} Q_N(n)
\end{equation}
Thus $A_N^+(T)$ is the occupied subset of $\T\times \R_+$ consisting 
of $T$ squares of area $1/N^2$.  We define  a {\em discrepancy function}
$D_{N,T}$ by 
\begin{equation} \label{discrepancy}
D_{N,T}(x,y) = N(1_{A_N(T)}- 1_{\{y\le T/N^2\}}).
\end{equation}
Figure \ref{symm diff} gives a closer look at the discrepancy 
between $A(T)$ and
the expected strip by distinguishing early and late sites relative to 
the time $T = yN^2$.  The figure 
depicts  the sign of $D_{N,T}$ in different 
colors.  $D_{N,T}$ takes on the values $\pm N$ and $0$.
$D_{N,T}>0$ means that $Q$ is early relative to the time $T$.  
$D_{N,T}<0 $ means that $Q$ is late.

The factor $N$ in the definition of $D_{N,T}$ is the appropriate 
normalization so that the limit exists in the sense of distributions
as $N\to \infty$.   Informally, our next theorem says that 
\begin{equation}
D_{N,T}(x,y) \to D(x) \d(y-y_0)
\end{equation}
in the sense of distributions with 
\[
D(x) \sim \  \sum_{k=1}^\infty 
\frac{a_k}{\sqrt{k}} \cos (2\pi kx) 
+
\frac{b_k}{\sqrt{k}} \sin (2\pi kx)
\]
and $a_k$ and $b_k$ independent, normally distributed random variables
with mean zero and variance $1$.    The random variable $D(x)$ is not 
defined for individual values of $x$.  For each $x$,  the variance, 
\[
\sum_{k=1}^\infty \frac1k(\cos^2(2\pi x) + \sin^2(2\pi x))  = \infty
\]
The precise statement of the theorem uses duality and involves
weight factors  that are merely asymptotic to $c/\sqrt k$.   

Let $H_0$ be the Sobolev space of functions $\eta$ on $\T\times \R_+$ 
satisfying $\eta(x,0) = 0$ and square norm equal to the Dirichlet integral,
\[
\|\eta\|_{H_0}^2 = \int_0^\infty \int_0^1 |\nabla \eta(x,y)|^2 \, dx dy 
\]
The restriction of $H_0$ to the circle $y=y_0$ is the Sobolev
space $H^{1/2}$.  Its dual is the space of $H^{-1/2}$ distributions
on $\T$ with dual norm given by
\[
\|f\|_{\{y_0\}} = \sup \{ \int_0^1 f(x) \eta(x,y_0) \, dx: \ \|\eta\|_{H_0} \le 1\}
\]

Fix an integer $K$, and consider test functions 
$\f\in C^\infty(\T \times \R)$ of the form
\[
\f(x,y) = \sum_{|k|\le K} \a_k(y) e^{2\pi i k x}
\]
Assume that for each $k$, $\a_k$ is supported in the annulus $0< c_1 \le |y| \le c_2$,
and the $\f$ is real-valued, i.~e.,  $\a_{-k} = \overline{\a_k}$.

\begin{theorem} \label{circle theorem} Let $T = \lfloor y_0 N^2 \rfloor$.  
Then as $N\to\infty$. 
\[ 
D_{N,T}(\f): = \int_{\T \times \R} D_{N,T}(x,y)\f(x,y)\, dxdy
\]
tends in law to a normally distributed random variable with mean zero and
variance
\[
S_{y_0}^2(\f) := \|\f(\cdot,y_0)\|_{\{y_0\}}^2 = \,  \sum_{0 <|k| \le K} m_k |\a_k(y_0)|^2
\]
with 
\[
m_k = \frac{1}{4\pi |k|} (1-e^{-4\pi |k| y_0})
\]
\end{theorem}
The messy term $e^{-4\pi |k|y_0}$ in the coefficients $m_k$ 
comes from starting the growth process at $y=0$.  If we started
at $y=-\infty$ it would disappear.  

In general, a gaussian random variable relative to a Sobolev space
has the form
\[
X = \sum_j a_j \f_j
\]
where $\f_j$ form an orthonormal basis for the Hilbert space
and $a_j$ are mean zero, variance $1$ independent random
variables.    Thus Theorem \ref{circle theorem} asserts that
$D_{N,T}$ tends to $D(x)\d(y-y_0)$ in
which $D$ is a (real-valued) gaussian random variable with mean value
zero associated to the Hilbert space of functions $g\in H_{y_0}^{1/2}(\T)$ with 
\[
\hat g(k) = \int_0^1 g(x) e^{-2\pi kx} \, dx
\]
$\hat g(0)= 0$,  $\hat g(-k)= \overline{\hat g(k)}$ 
\begin{equation}\label{half norm}
\|g\|_{H_{y_0}^{1/2}}^2  =  \sum_{k\neq0} |\hat g(k)|^2/m_k
\end{equation}

Roughly speaking, a discrepancy $|D(x)|$ of size one means
that the particles arrive late or early by a unit distance in the original lattice 
or distance $1/N$ in the continuum cylinder.   To illustrate this
we consider the example, with $T=N^2$, in which $A_N(T)$ occupies all the 
squares of $y\le 1-1/N$, 
exactly half of the $N$ squares in $1-1/N \le y \le  1$ and exactly
half of the $N$ squares in $  1 \le  y \le  1 + 1/N$.   Then 
$D_{N,T}(x,y)=N$ on each of the occupied squares of $1 < y < 1 + 1/N$ 
and $D_{N,T}(x,y) = -N$ on each of the unoccupied squares of 
$1 - 1/N < y < 1$.   In both cases the
integral of $|D_{N,T}|$ over the square is $N/N^2 = 1/N$ and there are $N$
such squares so the total is 
\[
\int_{\T\times \R} |D_{N,T}(x,y)| dx dy =  1
\]
Thus, in this example, the limit satisfies $|D(x)| = 1$ (half positive and half negative),

With the appropriate interpretation of the size of $D$ in mind, we can confirm heuristically the predictions of Meakin and Deutch as follows.   At scale
$N$, it's natural to truncate the series to $k\le N$, and say
\[
D(x) \sim  \   \sum_{k=1}^N
\frac{a_k}{\sqrt{k}} \cos (2\pi kx) 
+
\frac{b_k}{\sqrt{k}} \sin (2\pi kx)
\]
with independent unit variance gaussians as coefficients. 
For each fixed $x$, the variance of the right side is
\[
\sum_{k=1}^N \frac1k \approx \log N
\]
Thus the standard deviation of $D(x)$ is expected to be
on the order of $\sqrt{\log N}$.   On the other hand, we
can also predict the maximum fluctuation over all $x$.  
At scale $1/N$, we have $N$ different values of $x$
at which the discrepancy is represented by a random variable
of standard deviation  $\sqrt{\log N}$.   While these are not
independent, they are less and less correlated as the separation
gets larger.  Thus we expect the largest of $D(x)$ and the largest $-D(x)$
to be on the order of a factor $\sqrt{\log N}$ larger than a single standard
deviation, or $(\sqrt{\log N})^2 = \log N$.  This is the maximum bound 
demonstrated in Theorems \ref{maxfluc} and \ref{maxfluccyl}. 
The same heuristic reasoning applies in higher dimensions.  
The central limit theorems of  \cite{JLS11} in dimensions $d\ge 3$
yield a truncated variance of size $O(1)$ at typical boundary sites
consistent with the higher-dimensional numerical evidence of
Meakin and Deutch.   Moreover, the same reasoning as above predicts
that the maximum fluctuation in dimensions $d\ge 3$ is $O(\sqrt{\log T})$, where $T=N^d$ or $T=r^d$ is the number of particles.  This higher dimensional
estimate is proved in \cite{JLS2} and \cite{AG2}.  Very
recently in \cite{AG3}, Asselah and Gaudilli\`ere have confirmed that size
$\sqrt{\log T}$ fluctuations do occur.

\begin{figure}[htbp]
\begin{center}
\begin{tabular}{cc}
\includegraphics[width=.49\textwidth]{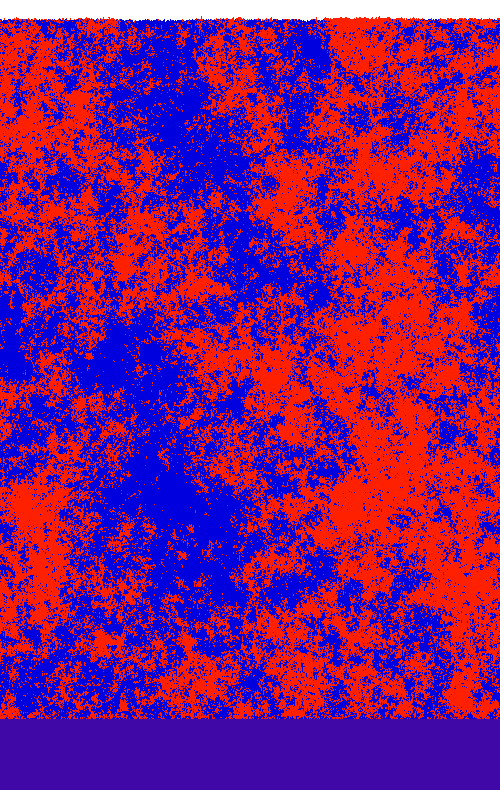}
&
\includegraphics[width=.49\textwidth]{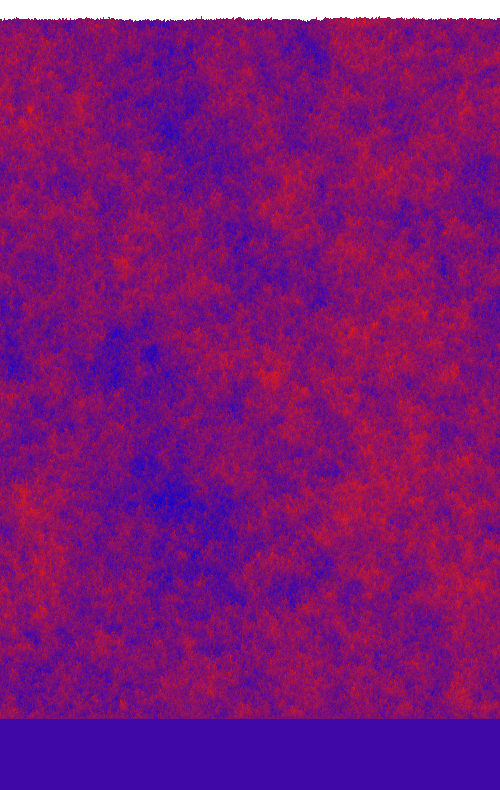}
\\
(a) & (b)
\end{tabular}
\caption{\label{whole symm diff}
Simulation of internal DLA on the cylinder with $N=500$. 
In (a) late points ($L_N>0$) and early points ($L_N<0$) 
are indicated in blue and red, respectively.  (b) The intensity of
the colors indicates the size of $L_N$.
}
\end{center}
\end{figure}

We will also analyze the fluctuations of the entire process as opposed
to what happens at a single time $T$.   In analogy with the discrepancy
function $D_{N,T}(x,y)$, we will define a rescaled
lateness function $L_N(x,y)$ that measures how early or
late the point $(x,y)$ is reached by the cluster.  

Before doing so, we will 
introduce a continuous time parameter $t$.
Let $T(t)$ be a standard Poisson random variable, with $T(0)=0$.   
$T(t)$ is an integer-valued process that produces in expectation
$t$ particles at time $t$.  For 
every $0 \le t_0 < t_1 < \cdots < t_m$, the variables
$T(t_j) - T(t_{j-1})$ are independent nonnegative integer-valued with
expectation $t_j- t_{j-1}$, respectively.  We also 
assume that $T(t)$ is independent of the internal DLA growth process,
and consider the process $A(T(t))$ depending on continuous 
time $t$.

For $n=(n_1,n_2)\in \Z_N\times \Z$, define
\begin{equation}
F(n) = \inf \{ t: n\in A(T(t)) \}
\end{equation}
and
\begin{equation}
L_N(x,y) = F(n)/N - yN,\quad (x,y) \in Q_N(n)
\end{equation}
For example, if $F(n) = (n_2+1)N$, then $n=(n_1,n_2)$ joins the cluster 
exactly one row late, in other words, 
by $N$ units of time $t$, which corresponds to a single row of width $N$ 
and height $1$ in $\Z_N\times \Z$ or a single row of width $1$ and height $1/N$ in
$\T\times \R$.   In that case, $L_N(n/N) = 1$.  Figure \ref{whole symm diff} depicts
simulations of $L_N$.



Informally, we say that $L_N$ tends as $N\to \infty$ to
the gaussian free field with Dirichlet boundary conditions on
$\T \times \R^+$, that is a gaussian random variable with
respect to the Hilbert space $H_0$.
The rigorous statement in dual form is as follows.  Let $\f(x,y)$ be
defined as above. 
\begin{theorem}  \label{gff}  As $N\to\infty$,
\[
L_N(\f) := \int_{\T \times \R} L_N(x,y)\f(x,y) \, dx dy
\]
tends in law to a mean zero gaussian random variable
with variance 
\[
S^2(\f): = \| \f \|^2_{H_0^*} =  
 \sum_{|k| \le K} \int_0^\infty \left| \a_k(y') e^{2\pi k(y-y')} dy'\right|^2 dy
\]
\end{theorem}
We will establish an estimate on the error in the central limit theorem 
of order $O(N^{-2/15})$, depending on the size of $\f$ and  the magnitude of $K$.

\section{IDLA and Hele-Shaw flow}

In this section, we give a heuristic description of the relationship between
 internal DLA and
the Hele-Shaw model.  This section contains no proofs, only formal
derivations.  The proof that the deterministic limit of internal DLA 
is Hele-Shaw flow, given in 2009 by Levine and Peres \cite{LP10}, 
proceeds via a discrete version of a classical obstacle problem.
We will discuss their work in slightly more detail in Section \ref{obstacle}.

Recall that in internal DLA from a single source in $\Z^2$
a particle takes a random walk from the origin 
in $A(T)$ until the first time it exits.  Then it stops
and augments the cluster to form $A(T+1)$.
It stops at sites $y$ at unit distance from $A(T)$ with probability $p_T(y)$, 
and this first exit probability is the discrete harmonic measure.  
In other words, it satisfies
\begin{equation} \label{discrete conservation}
v(0) =  \sum_{y\in \Z^2} v(y) p_T(y)
\end{equation}
for every function $v:\Z^2\to \R$ satisfying
\[
\ll v(x) = 0, \quad \mbox{for all} \ x\in A(T),
\]
where $\ll$ is the discrete Laplacian defined by 
\[
\ll f(x) = \frac14[f(x+ e_1)+ f(x-e_1) + f(x+e_2) + f(x-e_2)] - f(x)
\]
This suggests that the deterministic continuum
limit of the growth process is governed by harmonic measure.
Indeed, the continuum limit of the random walk is Brownian
motion, and, according to Kakutani's theorem, the hitting probability of Brownian 
motion starting from a point of a domain is the harmonic measure relative
to that point.  The continuum process in which a region grows proportionally
to its harmonic measure is known as Hele-Shaw flow.

Hele-Shaw flow describes the flow of fluid between two nearby
parallel plates.  The occupied region is essentially two-dimensional, so 
it is modeled by an open set $\O_t\subset \R^2$ at time 
$t$.  Given a domain $\O_0$ at time $t=0$, fluid
is pumped in at the origin so that the area 
grows at a uniform speed, $|\O_t| = t + |\O_0|$.  The
pressure $p(x,t)$ satisfies  $\D p(x,t) = -\d$ in $\O_t$ and $p(x,t)= 0$
on $\partial \O_t$.  The Hele-Shaw equation governing the growth
says that the normal velocity of the boundary of $\O_t$ 
is $|\nabla p|$.  Since $p$ is Green's function
for $\O_t$ with pole at the origin, Hele-Shaw's equation
can also be expressed as saying that the growth of the domain
is proportional to its harmonic measure.   
The correspondence with the discrete case is $\O_0 \leftrightarrow A(T_0)$,
$\O_t \leftrightarrow A(T_1)$, $|\O_0| = T_0/N^2$, $|\O_t| = T_1/N^2$,
and $t = (T_1-T_0)/N^2$.  

One way to solve the Hele-Shaw equation is to solve instead for
\begin{equation}\label{udef}
u(x,t) = \int_0^t p(x,s) ds.
\end{equation}
It's well known (c.~f.  \cite{GV})
that for each fixed $t$, $u$ solves an obstacle problem as follows.
Choose $\g(x,t)$ to be a function on 
$\R^2$ solving $\D \g= t\d + 1_{\O_0} - 1$.   Let $w$ solve the obstacle
problem
\[
w(x,t) = \inf \{ f : \D f \le 0, \ f\ge \g\}
\]
Although $w$ depends on the choice of $\g$, 
the set
\[
\O_t = \{x\in \R^2: w(x,t) > \g (x,t)\}
\]
and the function 
\[
u(x,t) = w(x,t)- \g(x,t) \ge 0
\]
are independent of the choice of $\g$.  On $\O_t$,
$ \D u =  -\D \g = 1- t\d -1_{\O_0}$, and on $\O_t^c$,
$u=0$.  In fact,  
\begin{equation}\label{ueq}
\D u = 1_{\O_t}  - 1_{\O_0} - t\d
\end{equation}
in all of $\R^2$.  

Conversely, starting from $u$, differentiate \eqref{ueq} with respect to $t$, 
to obtain 
\[
\frac{\partial }{\partial t} \D u(x,t) = V \s_t - \d
\]
where $V$ is the normal velocity of the boundary of $\O_t$
and $\s_t$ is the arc length measure of $\partial \O_t$.   
Define 
\begin{equation}\label{pdef}
p(x,t) = \frac{\partial }{\partial t} u(x,t)
\end{equation}
Then 
\[
\frac{\partial }{\partial t} \D u(x,t) = \D p(x,t) = -\d + |\nabla p| \s_t,
\]
and hence $p = (\partial/\partial t)u$ is the pressure for a Hele-Shaw fluid
cell with normal velocity $V = |\nabla p|$.  

The formulas above yield  conservation laws,
\begin{equation}\label{conservation}
v(0) = \frac{\partial}{\partial t} \int_{\O_t} v(x) \, dx
\end{equation}
for every harmonic function $v$.   We derive \eqref{conservation} in
integrated form by multiplying \eqref{ueq} by $v$ and integrating to obtain 
\[
0 = \int_{\R^2} (\D v) u \, dx = \int_{\R^2} v (\D u) \, dx  = \int_{\O_t} v \, dx - \int_{\O_0} v\, dx - tv(0)
\]
(One sees formally that the integration by parts has no boundary terms
because $u$ vanishes to second order on $\partial \O_t$ and is identically
zero outside.)  These formulas are also known as quadrature formulas \cite{GV}.

We have now come nearly full circle.  Let $\o_t$ be the harmonic measure
of $\O_t$ with respect to the origin, defined by the property
\[
v(0) = \int_{\partial \O_t} v(x) \o_t(dx)
\]
for every harmonic function in $\O_t$ with, say, continuous boundary
values.  Then $\o_t = |\nabla p| \s_t = V\s_t$, where $V$ is 
the normal velocity of $\partial \O_t$, and 
\[
v(0) = \frac{\partial}{\partial t} \int_{\O_t} v(x) \, dx = 
\int_{\partial \O_t} v(x) V \s_t(dx) = \int_{\partial \O_t} v(x) \o_t(dx)
\]
The discrete analogue is the equation we started with, \eqref{discrete conservation}.

For any fixed discrete harmonic function $v$, define
\[
M(T) = \sum_{n\in A(T)} v(n)
\]
If $v(0)=0$, then \eqref{discrete conservation} implies that
the conditional expectation of $M(T+1)$ given $A(T)$ is
\begin{equation}\label{mart1}
\ee(M(T+1)| A(T)) = \sum_{y\in \Z^2} v(y) p_T(y) = v(0) = 0
\end{equation}
In other words, $M$ is a martingale.  Martingales of this
type for various choices of $v$ are the main tools in the proofs
of theorems about fluctuations.  The martingale property
is an immediate consequence of the discrete version of
Kakutani's theorem.   The continuum theorems won't be necessary
to us;  they just help us to gain intuition.

Finally, we carry out a heuristic derivation that suggests
the form of the central limit theorems concerning fluctuations.  Suppose
that the boundary is given by a perturbation of the disk, in polar coordinates,
\[
r < R + \e f(\theta), \quad f(\theta) = \sum_{k\in \Z} \a_k e^{ik\theta}
\]
with $\a_{-k} = \overline{\a_k}$. We calculate the linearization
of Hele-Shaw flow for perturbations of the disk.  The 
Hadamard variational formula 
says that the (first order in 
$\e$) change in the gradient of Green's function is 
minus the radial derivative of the harmonic extension of $f$,
\[
-\left.\frac{\partial}{\partial r} \sum_{k\in \Z} \a_k (r/R)^ke^{ik\theta} \right|_{r=R}
= -\sum_{k\in \Z} \frac{k\a_k}{R} e^{ik\theta} 
\]
The minus sign is very important.  When $f(\theta) >0$,  the perturbation
is farther from the origin than the location $Re^{i\theta}$ on the circle
and the harmonic measure is smaller than average, and fewer
particles than average accumulate near $Re^{i\theta}$.   This
deterministic aspect of the process 
that keeps the shape close to circular. 

Next, we guess as to the stochastic ingredients of the evolution. 
We propose that the modes vary independently.  We expect that for some
constant $c>0$, $k> 0$, $t= \pi r^2$, 
\begin{equation}\label{gff-coeff}
d\a_k = -k \a_k \frac{dr}{r} + c \, dB_k(\rho) = -k \a_k d\rho + c\, dB_k, 
\quad (\rho = \log r)
\end{equation}
with independent white noise (derivative $dB_k$ of a Brownian motion $B_k$)
of equal amplitude in each mode.    The term $-k\a_k d\rho$ represents the
deterministic drift back towards the disk coming from the calculation above.

With $c=1$, this is the stochastic differential equation that yields the Gaussian Free Field.   In \cite{JLS} we find instead that 
the stochastic differential equation turns out to be
\[
d\a_k = -(k+1) \a_k d\rho + dB_k,\quad \rho = \log r
\]
The fact that $k$ is replaced by $k+1$ is related to the curvature
of the boundary.  The circumference circle of the circle increases
with $r$, so there there is room for more particles at the larger radius,
and the modes decrease slightly more than given in the rough
calculation above.  On the other hand, in the case of the cylinder, the 
circumference of the boundary circle of reference remains constant,
and we show in this paper that we get exactly the Gaussian Free Field.

\section{Proof of Theorem \ref{circle theorem}}


Note first that if $g(x) = e^{2\pi i k x}$, $k>0$  and 
\begin{equation*}
u(x,y) = \begin{cases}
g(x) \sinh(2\pi ky)/\sinh(2\pi ky_0), \quad & 0\le y \le y_0 \\
 g(x) e^{-2\pi k(y-y_0)}, \quad &  y_0 \le y < \infty
\end{cases}
\end{equation*}
Then  the restriction norm 
\[
\| g\|_{H^{1/2}} = \inf\{ \|v\|_{H_0}: v(x,y_0) = g(x)\}
\]
is achieved by the harmonic extension $u$.   This is proved by computing
\[
\int_0^\infty \int_0^1 |\nabla u|^2 \, dx dy = (4\pi k)/(1-e^{-4\pi ky_0}) = 1/m_k
\]
so that formula \eqref{half norm} holds.

Divide the outcomes of the cluster growth $A(T)$ into the three events.
Event 1, with probability at least $1-N^{-100}$ is
the event that the conclusion of Theorem \ref{maxfluccyl} holds, or,
put another way, $D_{N,T}$ is supported in the set 
\[
F= \{(x,y): |y- T/N^2| \le C (\log N)/N\}
\]
for all $T\le C_1N^2$.  Event 2, is the event that $D_{N,T} $ is supported
in $y\le C_2$ for all $T \le C_1N^2$, but Event 1 does not hold.   Thus
Event 2 has probability at most $N^{-100}$.  Event
3 is the complement of Events 1 and 2.  

To estimate the probability of Event 3, we recall from \cite{JLS} that
thin tentacles are rare events.  Lemma A of \cite{JLS} can be
stated in a nearly equivalent form as follows.

Denote $B(n)= \{m\in \Z_N\times \Z: n_2 - N/2 \le m_2  < n_2  + N/2\}$.   This
is a cylinder with about $N^2$ lattice sites. 
\begin{lemma}
\label{tentacles} (Thin tentacles)
There are positive absolute constants $C_0$, $b>0$,  and $c>0$ such that
for all $n \in \Z_N\times \Z$ with $n_2\ge N$, 
\begin{equation} \label{tentaclebound}
\pp \{ n\in A(T)  \  \mbox{and} \  ~ \# (A(T) \cap B(n) )  \leq b N^2 \}  \leq
C_0 e^{-c N^2/\log N}.
\end{equation}
\end{lemma}

This lemma implies that that for $C_2$ sufficiently large relative
to $C_1$, Event 3 has probability at most $O(e^{-cN^2/\log N})$.
Indeed, suppose there is $n\in A(T)$ such that $n_2 > C_2 N$ for
some $T\le C_1N^2$ .  Then since $\#A(T) = T$, and $A(T)$ is connected,
for at least one $n'\in A(T)$ 
with $n_2'\ge N$,  $\# B(n')\cap A(T) \le (2/C_2)C_1N^2$.   Thus
if $C_1/C_2 < b$, Lemma \ref{tentacles} applies to $B(n')$ and
Event 3 has probability at most $C_0C_2e^{-cN^2/\log N}$.

On Event 1 we will replace $\f$ by a harmonic function.
For $|k|\le K << N$, define $q(k,N)\ge 0$ by 
\begin{equation}\label{qdef}
1 -\cos(2\pi k/N) = \cosh(q/N) - 1
\end{equation}
It follows that 
\begin{equation}\label{qasympt}
q(k,N) = 2\pi |k| + O(1/N^2)
\end{equation}
Define for $n\in \Z/N\Z \times \Z$
\begin{equation}\label{psidef}
\p_0(n,T,N) = \sum_{0< |k| \le K} \a_k(T/N^2) e^{2\pi i n_1/N} e^{(q/N)(n_2-T/N)}
\end{equation}
The function $\p_0$ is discrete harmonic on the grid of lattice points with
spacing $1/N$ that equals an approximation to $\f - \a_0$ on the 
circle $\{(x,T/N^2): x\in \T\}$.

We claim that on Event 1,
\begin{equation}\label{discretize}
\int_{\T\times \R} D_{N,T}(x,y) \f(x,y)\, dxdy =
\frac1{N} \sum_{n\in A^+(T)} \p_0(n,T,N) + O(\log N/N)
\end{equation}
To prove this first note that
\begin{align*}
\int_{\T\times \R} D_{N,T}(x,y) & \f(x,y)\, dxdy = 
\int_{\T\times \R} D_{N,T}(x,y) (\f(x,y) - \a_0)\, dxdy  \\
&  = \int_{F} D_{N,T}(x,y) (\f(x,y) - \a_0)\, dxdy  \\
&  = N\int_{F} 1_{A_N(T)}  (\f(x,y) - \a_0)\, dxdy  
\end{align*}
Let $F_N = \{n\in \Z/N\Z\times \Z: |n_2 - T/N| \le C\log N\}$. 
Next, for $n\in F_N$ and $(x,y)\in Q_N(n)$,
\[
|\f(x,y) - \p_0(n,T,N)| \le C_3(\log N)/N
\]
where $C_3 = 10C\max |\nabla \f|$.   Without loss of
generality, $C \log N$ is an integer.  Therefore $F\cap A_N(T)$ is
a union of squares of side $1/N$ and we can match
every such square with its corner lattice point and replace
replace $\f -\a_0$ by $\p_0$.  Thus we obtain
\[
N\int_{F} 1_{A_N(T)}  (\f(x,y) - \a_0)\, dxdy  
= \frac1N \sum_{n\in F_N} 1_{A(T)} \p_0(n,T,N) + O((\log N)/N)
\]
Moreover,
\[
\frac1N \sum_{n\in F_N}1_{A(T)} \p_0(n,T,N) 
= 
\frac1N \sum_{n\in A^+(T)} \p_0(n,T,N) 
\]
This concludes the proof of \eqref{discretize}.  

Define
\[
M(s) = \frac1N\sum_{n\in A(T\wedge s)} \p_0(n,T,N)
\]
Then $M$ is a martingale.  Denote by
\[
Q = \sum_{s=1}^{T} \ee(|M(s)- M(s-1)|^2| A(s-1))
\]
the quadratic variation of the martingale, and denote
\[
S^2= \ee(Q), \quad B = \ee(|Q - S^2|^2), \quad A = \sum_1^T \ee|M(s) - M(s-1)|^4
\]
A theorem of Heyde and Brown \cite{HB}  gives a bound on the rate
of convergence in the martingale central limit theorem as follows.
There is an absolute constant
$C$ such that 
\[
\sup_{\l\in \R} |\pp(M(T)/S   \le \l) - \F(\l)| \le C \left( (A+B)/S^4\right)^{1/5}
\]
Note that 
\[
Q =  \frac1{N^2} \sum_{n\in A^+(T)} |\p_0(n,T,N)|^2 
\]
Define
\[
H(x,y) =  \sum_{0 < |k| \le K} \a_k(y_0) e^{2\pi i k x} e^{2\pi |k| (y-y_0)}
\]
On Event 1, $A^+_N(T)$ is up to a strip of unit width and height $C(\log N)/N$, 
equal to the the set $0 \le y \le y_0$.   Moreover, because $q = 2\pi |k| + O(1/N^2)$,
for $(x,y)\in Q_N(n)$,
\[
|H(x,y) - \p_0(n,T,N)|  = O(1/N)
\]
Thus on Event 1,
\[
Q =  \int_0^1 \int_0^{y_0} |H(x,y)|^2 \, dxdy + O(1/N)
\]
Furthermore,
\[
 \int_0^1 \int_0^{y_0} |H(x,y)|^2 \, dxdy  =  \sum_{0 < |k| \le K} m_k|\a_k(y_0)|^2 
\]
with 
\[
m_k =  \int_0^{y_0} e^{4\pi|k|(y-y_0)} \, dy = \frac1{4\pi |k|}(1-e^{-4\pi|k|y_0})
\]
Hence $|Q - S_{\{y_0\}}^2(\f)| \le 1/\sqrt{N}$ with probability $1-N^{-100}$.   
On Event 2, $A(T) \subset \{ n_2 \le C T/N\}$, so
that the factor $\disp e^{(q/N)(200N- T/N)} \le e^{200K}$ is bounded,
and  $Q = O(N^2)$.  Thus the expectation from Event 2 is at most 
$O(N^2 N^{-100}) = O(N^{-98})$.  Finally, on Event 3, the
worst case,  we still have the trivial estimate  $n\in A(T)\implies
n_2 \le T$.  Hence $e^{(q/N)(n_2- T/N)} \le e^{CN}$  for  constant $C$ depending
only on $K$, and $Q = O(e^{CN})$.  But Event 3 has probability of
order $e^{-cN^2/\log N}$, which is much smaller than exponential.
All together we have $\ee | Q  - S_{\{y_0\}}^2|^2 \le 1/N$.

On Event 1 or 2, $|M(s)-M(s-1)|^4 \le C/N^4$, 
This contributes to $B$ a sum of size  
$O(N^2 /N^4) =  O( 1/N^2)$.     On Event 3, the worst size
case is size $e^{CN}$ which is much smaller than $e^{-cN^2/\log N}$,
and hence negligible in the sum representing $B$. 

In all, $A+B = O(1/N)$ and we get the
bound $N^{-1/5}$ for the discrepancy of the distribution with the one for
the standard normal variable.  

\section{Proof of Theorem \ref{gff}}

We consider separately Events 1, 2, and 3 as in the proof of Theorem \ref{circle theorem}.    
\begin{lemma} \label{approx lateness} Denote
\[
L_N(\f)
 = \int_{\T \times \R}  L_N(x,y) \f(x,y) \, dx dy 
 \]
 On Event 1, with probability at least $1-N^{-100}$,
  \begin{align*}
L_N(\f)
&= 
\frac1{N^3} \int_0^\infty (t- T(t)) \a_0(t/N^2) dt \\
& \quad + \frac1{N^3} \int_0^\infty \sum_{n\in A^+(T(t))} \p_0(n,t,N) dt
+ O((\log N)^2/N)
\end{align*}
\end{lemma}

\noindent
Proof of Lemma \ref{approx lateness}.  
Rewrite $L_N$ as 
\[
L_N(x,y) = \int_0^\infty (1_{y\le T(t)/N^2} - 1_{A_N(T(t))} )dt
+ \int_0^\infty (1_{y\le t/N^2} - 1_{y\le T(t)/N^2}) dt
\]
It follows that
\[
\int_{\T \times \R}  L_N(x,y) \f(x,y) \, dx dy = U_1 + U_2 + U_3
\]
with 
\begin{align*}
U_1 &= 
\frac1{N} \int_0^\infty \int_0^1 \int_0^\infty (1_{\{y\le T(t)/N^2\}} -1_{A_N(T(t))})
(\f(x,y)-\a_0(y)) \, dt dx dy \\
U_2 &= \frac1{N} \int_0^\infty \int_0^1 \int_0^\infty (1_{\{y\le T(t)/N^2\}} -1_{A_N(T(t))})
\a_0(y) \, dt dx dy \\
U_3 &= \frac1{N} \int_0^\infty \int_0^\infty (1_{\{y \le t/N^2\}} - 1_{\{y\le T(t)/N^2\}} )
\a_0(y) \, dt dy
\end{align*}

Almost surely, $|T(t) - t| = O(t^{1/2} \log t)$.  Furthemore, on Event 1,
the integrand of $U_1$   is supported within
distance $O((\log N)/N$ of $y=t/N^2$. Therefore,  on the support,
\begin{equation}\label{harmonicreplace}
|\f(x,y)-\a_0(y) - \p_0(xN, yN, t, N)| = O((\log N)/N)
\end{equation}
Hence, when we replace $\f(x,y)-\a_0(y)$ by $\p_0(xN,yN,t,N)$,
the difference is dominated by 
\begin{align*}
\frac1{N} \int_{y_0N^2/2}^{2y_1N^2} \int_0^\infty  \int_0^1 
& |1_{\{y\le T(t)/N^2\}} -1_{A_N(T(t))}| ((\log N)/N ) \,  dx dy dt \\
& \le
\frac1{N} \int_{y_0N^2/2}^{2y_1N^2} ((\log N)/N)^2 \, dt 
= O((\log N)^2/N)
\end{align*}
Next, we claim that 
\begin{align*}
\frac1{N} &\int_0^\infty \int_0^\infty \int_0^1 
(1_{\{y\le T(t)/N^2\}} -1_{A_N(T(t))})
\p_0(xN,yN,t,N) \,  dx dy dt \\
& =
\frac1{N} \int_0^\infty \int_0^\infty \int_0^1 (1_{\{y\le h(t)\}}- 1_{A_N(T(t))})
\p_0(xN,yN,t,N) \,  dx dy dt  \\
& =
-\frac1{N^3} \int_0^\infty \sum_{n\in \Z_N\times \Z} (1_{\{n_2\le h(t)N\}} - 1_ {A^+(T(t)) })
\p_0(n, t,N) \,  dx dy dt  \\
& \qquad \qquad + O((\log N)/N) \\
& =
-\frac1{N^3} \int_0^\infty \sum_{n\in A^+(T(t)) }
\p_0(n, t,N) \,  dx dy dt  + O((\log N)/N)
\end{align*}
Indeed, the first equation is valid because 
\[
\int_0^1 \p_0(xN,yN,t,N) dx = 0
\]
In other words,
we may replace $1_{y\le T(t)/N^2}$ with $1_{y\le h(t)}$ for
any function  $h$.   In order to justify the error bound in
the second equation, choose $h(t) = \lfloor t/N\rfloor/N$.  
With probability $1-N^{-100}$ this, once again, confines the integration 
to a strip of  width $O((\log N)/N$.  
Replace 
$\p_0(xN,yN,t,N)$ on $Q_N(n)$ 
with the value at the corner $\p_0(n,t,N)$ to obtain a discrete sum.
(For convenience, we chose $h(t)$ so that the corresponding  discrete 
upper bound $n_2\le Nh(t) = \lfloor t/N\rfloor$ is an integer.)
Replacing $(xN,yN)$ with $n$ moves the point by  a distance
at most $1$ in each variable, which changes $\p_0$ by  $O(1/N)$. 
In the previous substitution, 
the difference in \eqref{harmonicreplace} was $O((\log N)/N)$,
so the integrated error here is smaller by the factor $1/\log N$.
Finally, in the last equation, replacing $h(t)$ by $0$ does not 
change the sum because for each $n_2$, 
\[
\sum_{n_1=1}^N \p_0(n_1,n_2,t,N) = 0
\]

Next we confirm that  with probability $1-N^{-100}$, 
\begin{equation}\label{E2bound}
U_2 = O((\log N)^2/N)
\end{equation}
Let
\[
f_N(t,y) = \int_0^1 (1_{\{y\le T(t)/N^2\}} - 1_{A_N(T(t))}) dx
\]
Since $A^+(T(t))$ consists of $T(t)$ squares of side length $1/N$,
\begin{equation}\label{fNcancel}
\int_0^\infty f_N(t,y) dy = 0
\end{equation}
On Event 1, and using the almost sure estimate $|T(t) - t| = O(t^{1/2}\log t)$, 
we have
\begin{equation}\label{fNsupport}
f_N(t,y) = 0 \quad \mbox{for}  \quad |y-t/N^2| \ge C(\log N)/N
\end{equation}
Consequently, since $|f_N(t,y)| \le 2$, 
\begin{equation}\label{fNbound}
\int_0^\infty |f(t,y)| dy = O((\log N)/N)
\end{equation}
Hence, on Event 1,
\begin{align*}
U_2 &= \int_0^\infty \int_0^\infty \frac1N f_N(t,y) \a_0(y) \,  dy dt \\
& = \int_{y_0N^2/2}^{2y_1N^2} \int_0^\infty \frac1N f_N(t,y) \a_0(y) \,  dy dt \\
& = \int_{y_0N^2/2}^{2y_1N^2} \int_0^\infty \frac1N f_N(t,y) (\a_0(y)- \a_0(t/N^2) 
\,  dy dt \\
& \le \int_{y_0N^2/2}^{2y_1N^2} ((\log N)/N^2) ((\log N)/N)  \, dt  = O((\log N)^2/N)
\end{align*}
The second equation the fact that $\a_0$ is supported in $y_0 \le y \le y_1$
and \eqref{fNsupport}.  
The third equation uses \eqref{fNcancel}.  The final inequality uses
\eqref{fNbound} and \eqref{fNsupport} to bound the difference in
values of $\a_0$.   This concludes the proof of \eqref{E2bound}.

Finally, we show that 
\begin{equation} \label{E3bound}
U_3 = \frac1{N^3} \int_0^\infty (t-T(t)) \a_0(t/N^2) \, dt + O((\log N)/N)
\end{equation}
Indeed, change variables to $r = N^2y$, to and define $R_1$ by 
\begin{align*}
U_3 &= \frac{1}{N^3}
\int_0^\infty \int_0^\infty (1_{\{ r\le t/N^2\}} - 1_{\{ r\le T(t)/N^2\}}) \a_0(r/N^2) \, dr dt\\
& = 
\int_0^\infty \int_0^\infty (1_{\{ r\le t/N^2\}} - 1_{\{ r\le T(t)/N^2\}}) \a_0(t/N^2) \, dr dt
+ R_1\\
& =\frac1{N^3} \int_0^\infty (t-T(t)) \a_0(t/N^2) \, dt + R_1
\end{align*}
Using $|T(t) - t| = O(t^{1/2} \log t)$, we see that 
\[
|\a_0(r/N^2) - \a_0(t/N^2)| = O((\log N)/N
\]
on the support of the integral representing the remainder term 
$R_1$.  This and the support properties of $\a_0$ yield
\[
|R_1| \le \frac1{N^3}\int_{y_0N^2/2}^{2y_1N^2} \int_{t-CN\log N}^{t+ CN\log N} ((\log N)/N) \, dr dt
 = O((\log N)^2/N)
 \]
 This concludes the proof of Lemma \ref{approx lateness} with an error of $O((\log N)^2/N)$.

For $0\le s <\infty$ define $M(s) = M_1 (s) + M_2(s)$ by
\begin{align*}
M_1(s) &= \frac1{N^3} \int_0^\infty  (t\wedge s - T(t\wedge s))\a_0(t/N^2) dt \\
M_2(s) &= \frac1{N^3} \int_0^\infty \sum_{n\in A^+(T(t\wedge s))} \p_0(n,t,N) dt
\end{align*}
The total quadratic variation of $M$ is 
\[
Q = \int_0^\infty Q(s) ds \quad \mbox{with} \quad Q(s) = \lim_{\e\to 0^+} \frac1\e
\ee( (M(s+\e)-M(s))^2| A(T(s)))
\]
Since $T(t)$ is independent of the process defining $A(T)$, $Q(s) = Q_1(s) + Q_2(s)$
the sum of the quadratic variation of $M_1$ and $M_2$ separately.
Because the $\a_k(y)$ are supported in $y\le y_1$,  
$M(s)$ is constant and $Q(s)=0$ for $s\ge y_1N^2$.

Let $M^*= \lim_{s\to\infty} M(s) = M(y_1N^2)$.   Then we have shown
that  with probability $1-N^{-100}$, 
\[
|L_N(\f) - M^*| \le  C(\log N)^2/N
\]

We will use the continuous analogue of the quantitative form of the martingale central limit theorem mentioned above, proved by Haeuser \cite{H}.   
\begin{lemma}\label{haeuser}  Let $S^2 = \ee Q$.  Then 
\[
\sup_{\l \in \R} |\pp(M^*/S \le \l) - \F(\l)| \le  \left[\ee(|Q-S^2|^2/S^4) + 
\ee \sum_{0\le s \le y_1N^2} |\Delta M(s)|^4/S^4\right]^{1/5}
\]
in which we define
\[
\sum_s |\Delta M(s)|^4 = \sum_j |M(s_j^+) - M(s_j^-)|^4
\]
\end{lemma}
This last sum is over the (almost surely finite number) of times $s_j$ in
$0\le s \le y_1N^2$ at which $M(s)$ is discontinuous.  
In fact, we will replace $S^2$ by $S^2(\f)$. 

\[
S^2(\f) = \sum_{|k| \le K} \int_0^\infty \left| \a_k(y') e^{2\pi k(y-y')} dy'\right|^2 dy
\]
In fact, these
$s_j$ are the times at which $T(s_j^+) - T(s_j^-) = 1$. 
\[
M_1(s_j^+) - M_1(s_j^-) = -\frac1{N^3} \int_{s_j}^\infty \a_0(t/N^2) \, dt
= -\frac1{N^3} \int_{s_j}^{y_1N^2} \a_0(t/N^2) \, dt
\]
\[
M_2(s_j^+) - M_2(s_j^-) = \frac1{N^3} \int_{s_j}^\infty \p_0(n_j, t,N)\, dt
= \frac1{N^3} \int_{s_j}^{y_1N^2} \p_0(n_j, t,N)\, dt
\]
with $\{ n_j\} = A(T(s_j^+)\backslash A(T(s_j^-))$.   On $\a_0= O(1)$ and
on Events 1 or 2, $p_0(n_j,t,N)$ is bounded so that $|\Delta M(s_j)| = O(1/N)$.
On Event 3, $|\Delta M(s_j)| = O(e^{CN})$, which is much smaller
than the probability $e^{-cN^2/\log N}$ of Event 3.  
There are are almost surely $O(N^2)$ jumps $s_j$.   Therefore,
\[
\ee \sum_{0\le s \le y_1N^2} |\Delta M(s)|^4/S^4 = O(1/N^2)
\]
We will show below that
\begin{equation}\label{quadvar4norm}
\ee(|Q-S^2(\f)|^2) = O(N^{-2/3})
\end{equation}
Once we have proved this, the proof of Theorem \ref{gff} is nearly
complete.  Lemma \ref{haeuser}  says
that $M^*$ has the same distribution as a gaussian with
variance $S^2=\ee Q$ up to $O(N^{-2/15})$.   Moreover, \eqref{quadvar4norm}
implies that $|S^2 - S^2(\f)| = O(N^{-1/3})$.   
According to Lemma \ref{approx lateness},
$M^*$ differs from $L_N(\f)$ by at most $O((\log N)^2/N)$ with probability $O(1- N^{-100})$.

It remains to prove \eqref{quadvar4norm}.  
Since $\ee((\e - (T(s+\e)-T(s)))^2) = \e$, 
\[
Q_1(s) = \left(\frac1{N^3} \int_s^\infty \a_0(t/N^2) dt\right)^2
\]
and
\[
Q_1 = \int_0^\infty Q_1(s) ds = \int_0^\infty \left(\int_y^\infty \a_0(y')dy'\right)^2 dy
\]

\begin{lemma} \label{quadvarlate} With probability $1-N^{-100}$, 
\[
Q_2 = \sum_{0 < |k| \le K} \int_0^\infty \left| \a_k(y') e^{2\pi k(y-y')} dy'\right|^2 dy
+ O(N^{-1/4})
\]
\end{lemma}
Proof.  
\begin{equation}\label{Q2}
Q_2(s) = \sum_{n\in \Z_N\times \Z} 
\left| \frac1{N^3} \int_s^\infty 
\p_0(n,t,N)\, dt \right|^2 p_s(n)
\end{equation}
where $p_s(n)$ is the probability that the random walk starting at
$(n_1,0)$,  $n_1= 1,\dots, N$ with equal probability, exits $A(T(s))$ 
for the first time at the site $n$.   Thus $p_s(n)$ is nonzero only
on the boundary of $A(T(s))$, that is at sites at distance exactly $1$ from 
 $A(T(s))$.    
 With probability $1-N^{-100}$, $p_s(n) >0$ 
 implies $|n_2 - s/N| \le C \log N$.  In other words, the 
 boundary of $A(T(s))$ is nearly a horizontal line.   One
 can therefore deduce from 
 barrier estimates (discrete harmonic majorants, not
 written down explicitly here)  
 that the distribution of $p_s(n)$ is approximately uniform
in the $n_1$ variable in the following sense.  
For any $a>0$, define
\[
R(a) = \{n: 1 \le n_1 \le aN, \ |n_2 - s/N| \le C \log N\}
\]
Then 
\[
\sum_{n\in R(a)} p_s(n) =  a + O(N^{-1/3})
\]
(Sharper bounds are also valid; we have not attempted to optimize
the power.)

Put the sites for which $p_s(n(j))>0$ in order according to their
position horizontally, $1\le n_1(1) \le n_1(2) \le \cdots $ and
consider disjoint intervals $I_j$ of $0\le x \le 1$, so
that the right endpoint of $I_j$ is the left endpoint of $I_{j+1}$
and the length $|I_j| = p_s(n(j))$.  Then for all $x\in I_j$,
\[
|\p_0(n(j),t,N)-  p_0(xN,s/N, t, N)| = O(N^{-1/3})
\]
Thus, if we replace the sum of $\p_0$ on the lattice with
the weighting $p_s(n)$ by the integral $dx$, we find that
 \[
Q_2(s) = \int_0^1 \left| \frac1{N^3} \int_s^{y_1N^2}
\p_0(xN,t,N)\, dt \right|^2 \, dx + O(N^{-2} N^{-1/3}) 
\]
The worst error comes from cross terms in the integrand
of the square with one factor of size $O(N^{-1/3})$ and the other of unit 
size.  This yields errors which are a factor $O(N^{-1/3})$ smaller
than the main term. The main term $Q_2(s)$ is of size $O(N^{-2})$
as one can see from the fact that 
the expression inside the $|\cdot|$ sign is an integral in $t$ 
over an interval of length of order $N^2$ of a (roughly) unit
sized integrand divided by $N^3$, thus of size $1/N$.  This
is squared and summed over the probability measure $p_s(n)$.

Finally, integrating $Q_2(s)$ over $0\le s \le y_1N^2$, and changing variables,
one finds the expression in Lemma \ref{quadvarlate}, 
with an error  that is a factor $O(N^{-1/3})$ smaller.

\section{Obstacle problems}\label{obstacle}

\begin{figure}
\begin{center}
\includegraphics[height=.3\textheight]
{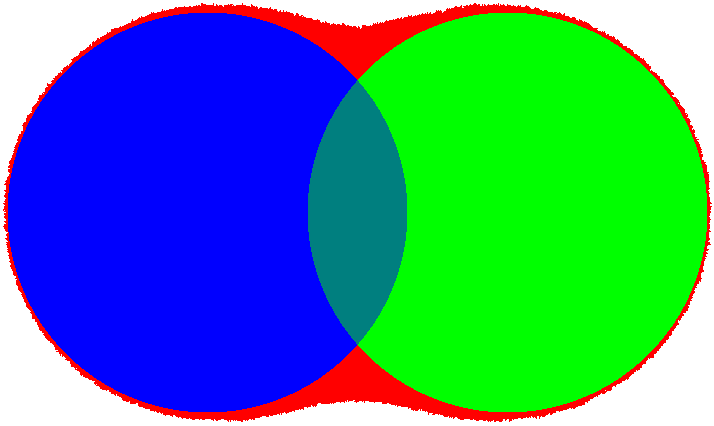}
\\ \vspace{3mm}
\caption{\label{smash sum}
Smash sum of two overlapping disks, $1_A + 1_B$. 
}
\end{center}
\end{figure}

Levine and Peres had an entirely different motivation.
In 1991 Diaconis and Fulton \cite{DF} defined the notion of {\em smash sum}.
Consider two open subsets $A$ and $B$ of the plane and define a function $\mu_N$ 
on the lattice $\Z^2$ by 
\[
\mu_N(n) = 1_A(n/N) + 1_B(n/N)
\]
The function $\mu_N$ represents an initial collection of particles, two
at each site of $A\cap B$ in a rescaled lattice of mesh size $1/N$ and
one particle at each of the rest of the sites in $A\cup B$.   Choose
any site with  more than one particle, and move one of them to any 
of the four adjacent sites with equal probability.    The order in which 
the particles move is irrelevant because they are interchangeable.  
A site can be occupied by many particles at the same time,
and the process is said to stop the first time each occupied site 
has exactly one particle.   Denote the final distribution by $\nu_N$.
The theorems  of Levine and Peres characterized the deterministic
limit of $\nu_N$, the set $C$ such that
\[
\lim_{N\to \infty} \frac{1}{N^2} \sum_{n\in \Z^2} \f(n/N) \nu_N(n)   = 
\int_{\R^2} \f(x) 1_{C}(x) dx
\]
Figure \ref{smash sum} gives an instance of the smash sum of two disks.   

\begin{figure}
\begin{center}
\includegraphics[height=.3\textheight]
{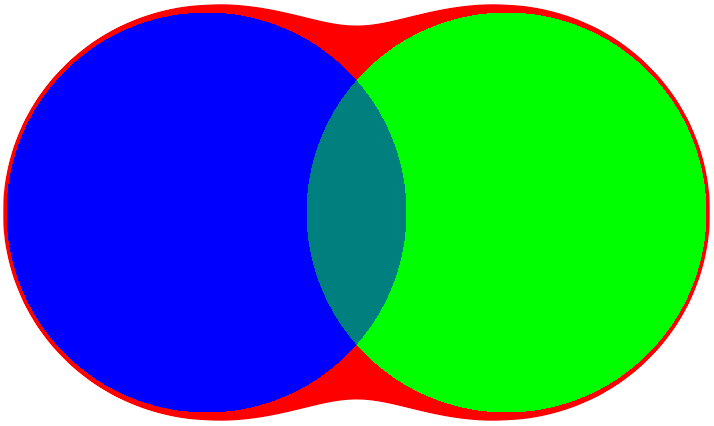}
\\ \vspace{3mm}
\caption{\label{div sandpile}
Divisible sandpile sum of two overlapping disks $1_A+ 1_B$.
}
\end{center}
\end{figure}

We will now give a few more of the ideas behind the
work of Levine and Peres \cite{LP10}.  The idea is to show
that the continuum limit solves the obstacle problem.   

The first step of their proof is
to study a deterministic process on the lattice, known as the divisible
sandpile.  (The analogous continuum process is balayage or sweeping out 
\cite{GV}.)   Consider any nonnegative function 
$\mu: \Z^2 \to \R_+$, representing possibly fractional quantities of particles 
at each lattice site or else the height
of a sandpile.   Pick a site $x$ at which  $\mu(x)>1$, after a single step
in the process, a toppling of the site $x$, the new heights are $1$ at $x$ 
and $(\mu(x)-1)/4$ extra at each of the four adjacent sites.  This
process continues until every site has at most height $1$.   
One can show that the process stops in a finite number of steps
for finitely supported $\mu$.   Moreover, the final height function, 
denoted $\nu(x)$, is independent 
of the order in which the toppling occurs.    The divisible
sandpile reflects the average behavior of the random walk.  
Figure \ref{div sandpile} depicts starting from the same $\mu_N$ as in
the preceding picture, Figure \ref{smash sum}.  The resemblance is already striking
for $N \approx 200$.   

Levine and Peres show that the random process is well approximated by the deterministic sandpile process by means of an auxiliary function they
call the odometer function $u(x)$.  The function $u(x)$ records the (fractional) number of particles donated by the site $x$ in the course of
the deterministic process.   The word odometer reminds
us that $u(x)$ does not represent a net loss
of particles, but rather the total quantity donated without subtracting
the number received.    It is not hard to check that $u$ solves the
discrete Laplace equation
\[
\ll u(x) = \nu(x) -\mu(x) 
\]
Moreover, $u$ can be obtained from the solution to the 
discrete obstacle problem as follows.  

\begin{lemma}\label{levine-peres}  (Levine-Peres) 
Fix any $\g(x)$ satisfying $\ll \g = \mu -1$.    Let $w$ solve the obstacle problem
\[
w(x) = \inf\{f(x): \ll f \le 0, \quad f\ge \g\}.
\] 
Let $u$ be the odometer function starting from $\mu$.
Then 
\[
u = w-\g.
\]
\end{lemma}
The fact that the discrete function $u$ tends to its continuum counterpart
depends ultimately on estimating the difference between the fundamental
solution of $\ll$ and $\D$.

\section{Further remarks and questions}

Thereom \ref{maxfluc} is proved in \cite{JLS} in dimension 2 and in
\cite{JLS2} in dimensions $d\ge 3$.  The proof uses martingales associated
to the discrete analogue of Green's function with a pole at a point near the putative boundary,  either inside or outside.   
These martingales take values 
larger than the expected value if there are extra points 
in the cluster near the pole, and the martingale is smaller than its expected
value if there are fewer then the typical number of occupied sites
in the cluster near the pole.   The central limit theorem is 
inadequate to the task of estimating large deviations of the martingale.  
Instead one uses the parametrization of the martingale by Brownian motion.  
The lemma concerning thin tentacles, Lemma \ref{tentacles}, is crucial as
well.   The proof also involves an iteration of successively better
estimates on the inner and outer deviations of the shape.
The proof of Theorem \ref{maxfluccyl} is analogous to that of Theorem
\ref{maxfluc}, just as Theorems \ref{circle theorem} and \ref{gff} are similar
to the corresponding theorems for $\Z^2$.

In this paper, we chose to treat the case of the cylinder so as
to identify a case in which the fluctuations are described exactly
by the Gaussian Free Field.  In \cite{JLS11} we carried out the
case of the disk.    One difference with the case of the cylinder is
that it's somewhat harder to construct suitable discrete harmonic
functions approximating $z^k$.  Furthermore, the estimates require
variants for averages with respect to discrete harmonic
functions of van der Corput's theorem counting lattice points in
disks.   We have not yet carried out the case $d=3$, although we
believe it follows from very similar methods.  The technical difficulty
is that it requires variants of theorems stronger theorems than van der 
Corput's concerning  the number of lattice points in a ball in $3$-space, 
along the lines of improvements due to Vinogradov (see \cite{IKKN}).

Whereas the square Dirichlet norm is 
\begin{align*}
\int_{\R^2} |\nabla f(x,y)|^2 \, dx dy 
&= 
\int_0^{2\pi} \int_0^\infty (|r\partial_r f|^2 + |\partial_\theta f|^2)\frac{dr}{r} d\theta  \\
&= 
2\pi \sum_{k\in\Z}
\int_0^\infty (|r\partial_r f_k|^2 + |kf_k|^2)\frac{dr}{r} 
\end{align*}
in which
\[
f(x,y)= \sum_{k\in\Z} f_k( r) e^{ik\theta}
\]
The square of the norm of the gaussian random field representing
fluctuations from a source at the origin in $\Z^2$ is 
\[
2\pi \sum_{k\in\Z}
\int_0^\infty (|r\partial_r f_k|^2 + |(|k|+1)f_k|^2)\frac{dr}{r} 
\]
In general, we expect that the random field will reflect the curvature
of the deterministic region.  But even in this simple case,
the norm is expressed in terms of non-local (pseudo-differential) 
operators.

The expression for the norm in Theorem \ref{circle theorem} 
at distance $y_0$, starting from the (exactly straight) boundary of 
$y\le 0$ involves  the  factor $(1-e^{-4\pi|k|y_0})$.   Thus in
some average sense, the influence of deterministic behavior at 
$y=0$  attenuates at an exponential rate at $y=y_0$.  It would be 
nice to understand this better.  At the same time one can
ask about the mixing time, that is, given a known boundary
at one time, how long do we need to wait before that configuration
is mostly forgotten?

One can also ask questions about random walks other than
the standard one.   The first author\footnote{The author thanks Pavel Etingof
for suggesting this problem.}  
supervised work on this subject in the summer of 2008 
by a high school student, Max Rabinovich \cite{R}.  He adapted 
the methods of Levine and 
Peres to the hexagonal lattice.  His key observation is
that the same methods work, provided one can approximate
the discrete fundamental solution by the analogue of the Newtonian 
potential.   At first it appears that the estimates need to
be good to second order at infinity which they are not for
the hexagonal lattice.  But on closer inspection, what is
required are error estimates for the difference of fundamental solutions,
as compared to the gradient of the Newtonian potential.  The
gradient is of order $1/r^{d-1}$ as $r\to\infty$ and the error
term is one order better, $1/r^d$, which is two orders better than $1/r^{d-2}$
as required.  Rabinovich's  theorem applies to all random walks on $\Z^d$ 
given by a finitely supported probability measure $p$ on $\Z^d$ 
such that the random walk moves from $x$ to $x+\alpha$ with probability
$p(\alpha)$ and 
\[
\sum_{\a \in \Z^d} p(\a)\a  = 0
\]
This condition means that the random walk has no drift.

It remained to consider walks with drift.    James Propp 
proposed the specific example of  a walk that moves 
East or North, each with probability $1/2$. 
If the source is the origin, this fills a cluster in the first quadrant.
If there are $T$ particles, it is natural to rescale by parabolic scaling 
$u=x+y$, and $v=x-y$ are replaced by $U= (x+y)/T^{2/3}$,  
$V = (x-y)/T^{1/3}$.  Then in parallel with the work of Levine and
Peres, one expects the cluster to be associated with an
obstacle problem based on parabolic operators in the $(U,V)$ 
variables as treated by Caffarelli, Petrosyan and
Shahgholian \cite{CPS}.  

After this lecture, Cyrille Lucas \cite{Lucas}, carried out this program and 
in the process established the existence
of a so-called heat ball.   Take the limiting (and indeed simplest case)
of the Hele-Shaw flow in which $\O_0$ shrinks to a point.  Then $\O_t$ is
the Euclidean ball of volume $t$.  The conservation law \eqref{conservation}
in integrated form can be written
\[
v(0) = \frac{1}{\mbox{vol} \, B} \int_{B} v(x) \, dx
\]
for any harmonic function $v$.  Of course, this is just the
well known mean value property for harmonic functions.
The domain analogous to the ball for the heat operator
is a set  $H\subset \{(x,t): t\ge 0\}$ 
of area $1$ such that 
\[
v(0,0) = \int_H u(x,t) \, dx dt
\]
where $v$ satisfies the adjoint or backwards heat equation
\[
[(\partial/\partial x)^2 + (\partial/\partial t)]v(x,t) = 0
\]
Evidently, any parabolic dilation $H_R = \{(Rx,R^2t): (x,t) \in H\}$ satisfies
of $H$ satisfies 
\[
v(0,0) = \frac{1}{R^3} \int_{H_R} v(x,t) \, dx dt
\]
Many weighted averages of $v$ produce $v(0,0)$.
This one is interesting because the weight is constant, proportional
to Lebesgue measure.

As mentioned in the lecture, the question that remains open
is the regularity of the boundary of $H$.   The discrete
construction of the divisible sandpile gives an approximation
to the continuum set $H$.  The theorems of
\cite{CPS} give a criterion involving approximations to $H$.
Their criterion would imply that the boundary of $H$
is smooth if there were a practical bound on the constants involved.
It would also be interesting to show that $H$ is convex, which
looks rather obvious from the sandpile approximation.  A typical
approach would be to realize $H$ as the level set of a log
concave function.  However, the
odometer function $u$ associated with $H$ is not log concave.
On the other hand, we have numerical evidence
that $u/h$ is log concave, in which $h$ is the standard fundamental
solution, $t^{-1/2}e^{-x^2/4t}$.    This would imply that $H = \{ u/h>0\}$ is convex.
It is not hard to show in the Propp example that the discrete analogue of $u/h$ 
is log concave in the $x$ direction.  Unfortunately, one does find numerically
a very few sites near the boundary at which log concavity fails
slightly in the $t$ direction.  So at least in the Propp 
example, it's hard to see how a combinatorial proof could succeed.

\end{document}